\documentclass{article}

\usepackage{amssymb,latexsym,amsmath}
\usepackage[pdftex]{graphicx}

\usepackage{hyperref}

\begin{document}

\pagestyle{plain}

\newtheorem{theorem}{Theorem}

\newtheorem{proposition}[theorem]{Proposition}

\newtheorem{lema}[theorem]{Lemma}

\newtheorem{corollary}[theorem]{Corollary}

\newtheorem{definition}[theorem]{Definition}

\newtheorem{remark}[theorem]{Remark}

\newtheorem{exempl}{Example}[section]

\newenvironment{exemplu}{\begin{exempl}  \em}{\hfill $\square$

\end{exempl}}  \vspace{.5cm}

\renewcommand{\contentsname}{ }

\title{Origin of emergent algebras}

\author{Marius Buliga \\ 
\\
Institute of Mathematics, Romanian Academy \\
P.O. BOX 1-764, RO 014700\\
Bucure\c sti, Romania\\
{\footnotesize Marius.Buliga@imar.ro}}  \vspace{.5cm}

\date{This version: 04.04.2013}

\maketitle

\begin{abstract}
This is an expository article concerning the last section "Why is the tangent space a group?" (section 8) of the  article by A. Bella\"{\i}che, The tangent space in sub-riemannian geometry, from the viewpoint of emergent algebras.
\end{abstract}

\section{Introduction}
\label{secone}

In the last section "Why is the tangent space a group?" (section 8) of the  article by A. Bella\"{\i}che, The tangent space in sub-riemannian geometry, the author explains a very interesting story, where names of Gromov and Connes appear, which is the first place, to my knowledge, where the idea of \href{http://arxiv.org/abs/0907.1520}{emergent algebras} appear.

Let me give the relevant passages.

[p. 73] "Why is the tangent space a group at regular points? [...] I have been puzzled by this question. Drawing a Lie algebra from the bracket structure of some $  X_{i}$'s did not seem to me the appropriate answer. I remember having, at last, asked M. Gromov about it (1982). The answer came under the form of a little apologue:

Take a map $  f: \mathbb{R}^{n} \rightarrow \mathbb{R}^{n}$. Define its differential as

$$ \mbox{(79)} \quad \quad  D_{x} f(u) = \lim_{\varepsilon \rightarrow 0} \varepsilon^{-1} \left[ f(x+\varepsilon u) - f(x) \right]$$

provided convergence holds. Then $  D_{x}f$ is certainly homogeneous:

$$  D_{x}f(\lambda u) = \lambda D_{x}f(u)$$

but it need not satisfy the additivity condition

$$  D_{x}f(u+v) = D_{x}f(u) + D_{x}f(v)$$

[...] However, if the convergence in (79)  is uniform on some neghbourhood of $  (x,0)$  [...]  then $  D_{x}f$ is additive, hence linear. So, uniformity was the key. The tangent space at $  p$ is a limit, in the [Gromov-]Hausdorff sense, of pointed spaces [...] It certainly is a homogeneous space -- in the sense of a metric space having a 1-parameter group of dilations. But when the convergence is uniform with respect to $  p$, which is the case near regular points, in addition, it is a group.

Before giving the proof, I want to tell of another, later, hint, coming from the work of A. Connes. He has made significant use of the following observation: The tangent bundle $  TM$ to a differentiable manifold $  M$ is, like $  M \times M$, a groupoid. [...] In fact TM is simply a union of groups. In [8], II.5, it is stated that its structure may be derived from that of $  M \times M$ by blowing up the diagonal in $  M \times M$. This suggests that, putting metrics back into the picture, one should have

$$\mbox{(83)} \quad \quad        TM = \lim_{\varepsilon \rightarrow 0} \varepsilon^{-1} (M \times M)$$

[...] in some sense to be made precise.

There is still one question. Since the differentiable structure of our manifold is the same as in Connes' picture, why do we not get the same abelian group structure? One can answer: The differentiable structure is strongly connected to (the equivalence class of) Riemannian metrics; differentiable maps are locally Lipschitz, and Lipschitz maps are almost everywhere differentiable. There is no such connection between differentiable maps and the metric when it is sub-riemannian. Put in another way, differentiable maps have good local commutation properties with ordinary dilations, but not with sub-riemannian dilations $  \delta_{\lambda}$.

{\em So, one should not be abused by (83) and think that the algebraic structure of $  T_{p}M$ stems from the absolutely trivial structure of $  M \times M$! It is concealed in dilations, as we shall now prove.}"

\section{Details and comments of Bella\"{\i}che's proposal}
\label{sectwo}

In this section we shall see how Bella\"{\i}che proposes to extract the algebraic structure of the metric  tangent space $  T_{p}M$ at a point $  p \in M$, where $  M$ is a regular sub-riemannian manifold. The metric tangent space is defined up to arbitrary isometries fixing one point, as the limit in the Gromov-Hausdorff topology over isometry classes of compact pointed metric spaces

$$  [T_{p} M, d^{p}, p] = \lim_{\varepsilon \rightarrow 0} [\bar{B}(p, \varepsilon), \frac{1}{\varepsilon} d, p]$$

where $  [X, d, p]$ is the isometry class of the compact  metric space $  (X,d)$ with a marked point $  p \in X$. (Bella\"{\i}che's notation is less precise but his previous explanations clarify that his relations (83), (84) are meaning exactly what I have written above).

A very important point is that moreover, this convergence is uniform with respect to the point $  p \in M$. According to Gromov's hint mentioned  by  Bella\"{\i}che, this is the central point of the matter. By using this and the structure of the trivial pair groupoid $  M \times M$, Bella\"{\i}che proposes to recover the Carnot group algebraic structure of $  T_{p}M$.

From this point on I shall pass to a personal interpretation of the  section 8.2 "A purely metric derivation of the group structure in $  T_{p}M$ for regular $  p$" of Bella\"{\i}che article. (We don't have to worry about "regular" points because I already supposed that the manifold is "regular", although Bella\"{\i}che's results are more general, in the sense that they apply also to sub-riemannian manifolds which are not regular, like the Grushin plane.)

In order to exploit the limit in the sense of Gromov-Hausdorff, he needs first an embodiment of the abstract isometry classes of pointed metric spaces. More precisely, for any $  \varepsilon  >  0$ (but sufficiently small), he uses a function denoted by $  \phi_{x}$, which he states that it is defined on $  T_{x} M$ with values in $  M$. But doing so would be contradictory with the goal of constructing the tangent space from the structure of the trivial pair groupoid and dilations. For the moment there is no intrinsic meaning of $  T_{x} M$, although there is one from differential geometry, which we are not allowed to use, because it is not intrinsic to the problem.  Nevertheless, Bella\"{\i}che already has the functions $  \phi_{x}$, by way of his lengthy proof (but up to date the best proof) of the existence of adapted coordinates. For a detailed discussion see my article "Dilatation structures in sub-riemannian geometry" \href{http://arxiv.org/abs/0708.4298}{arXiv:0708.4298}.

Moreover, later he mentions "dilations", but which ones? The natural dilations he has from the vector space structure of the tangent space in the usual differential geometric sense? This would have no meaning, when compared to his assertion that the structure of a Carnot group of the metric tangent space is concealed in dilations.  The correct choice is again to use his adapted coordinate systems and use intrinsic dilations.  In fewer words, what Bella\"{\i}che probably means is that his functions $  \phi_{x}$ are also decorated with the scale  parameter $  \varepsilon > 0$, so they should deserve the better notation $  \phi_{\varepsilon}^{x}$,  and that these functions behave like dilations.

A natural alternative to Bella\"{\i}che's proposal would be to use an embodiment of the isometry class $  [\bar{B}(x, \varepsilon), \frac{1}{\varepsilon} d, x]$ on the space $  M$, instead of the differential geometric tangent space $  T_{x}M$.  With this choice, what Bella\"{\i}che is saying is that we should consider dilation like functions $  \delta^{x}_{\varepsilon}$ defined locally from $  M$ to $  M$ such that:
\begin{enumerate}
	\item[-]they preserve the point $  x$ (which will become the "$  0$" of the metric tangent space): $  \delta^{x}_{\varepsilon} x = x$ 
	\item[-]they form a one-parameter group with respect to the scale: $  \delta^{x}_{\varepsilon} \delta^{x}_{\mu} y = \delta^{x}_{\varepsilon \mu} y$ and $  \delta^{x}_{1} y = y$, 
	\item[-]for any $  y, z$ at a finite distance from $  x$ (measured with the sub-riemannian distance $  d$, more specifically such that  $  d(x,y), d(x,z) \leq 1$) we have 
\end{enumerate}
$  d^{x}(y,z) = \frac{1}{\varepsilon} d( \delta^{x}_{\varepsilon} y, \delta^{x}_{\varepsilon}z) + O(\varepsilon)$

where $  O (\varepsilon)$ is uniform w.r.t. (does not depend on) $  x, y , z$ in compact sets.

Moreover, we have to keep in mind that the "dilation"  $  \delta^{x}_{\varepsilon}$ is defined only locally, so we have to avoid to go far from $  x$, for example we have to avoid to apply the dilation for $  \varepsilon$ very big to points at finite distance from $  x$.

Again, the main thing to keep in mind is the uniformity assumption. The choice of the embodiment provided by "dilations" is not essential, we may take them otherwise as we please, with the condition that at the limit $  \varepsilon \rightarrow 0$ certain combinations of dilations converge uniformly. This idea suggested by Bella\"{\i}che reflects the hint by Gromov.  In fact this is what is left from the idea of a manifold in the realm of sub-riemannian geometry  (because adapted coordinates cannot be used for building manifold structures, due to the fact that "local" and "infinitesimal" are not the same in sub-riemannian geometry, a thing rather easy to misunderstand until you get used to it).

Let me come back to Bella\"{\i}che reasoning, in the setting I just explained. His purpose is to construct the operation in the tangent space, i.e. the addition of vectors. Only that the addition has to recover the structure of a Carnot group, as proven by Bella\"{\i}che. This means that the addition is not a commutative, but a noncommutative  nilpotent operation.

OK, so we have the base point $  x \in M$ and two near points $  y$ and $  z$, which are fixed. The problem is how to construct an intrinsic addition of $  y$ and $  z$ with respect to $  x$. Let us denote by $  y +_{x} z$ the result we are seeking. (The link with the trivial pair groupoid is that we want to define an operation which takes $  (x,y)$ and $  (x,z)$ as input and spills $  (x, y+_{x} z)$ as output.)

The relevant figure is the following one, which is an improved version of the Figure 5, page 76 of 
Bella\"{\i}che paper.

\vspace{.5cm}   \centerline{\includegraphics[width=120mm]{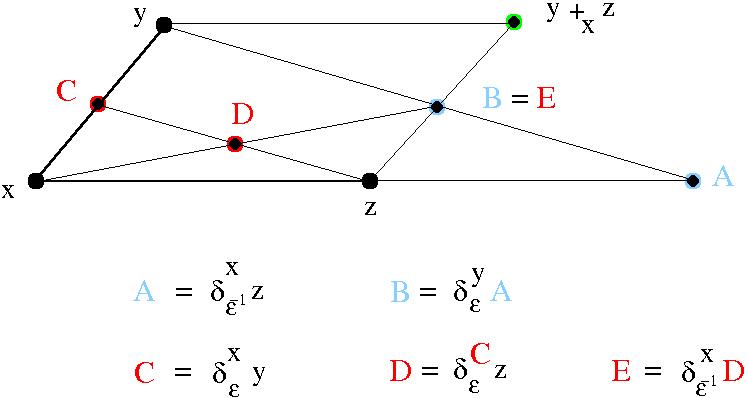}}  \vspace{.5cm}

Bella\"{\i}che's recipe has to do with the points in blue. He says that first we have to go far from $  x$, by dilating the point $  z$ w.r.t. the point $  x$, with the coefficient $  \varepsilon^{-1}$. Here $  \varepsilon$ is considered to be small (it will go to $  0$), therefore $  \varepsilon^{-1}$ is big.  The result is the blue point $  A$. Then, we dilate (or rather contract) the point $  A$  by the coefficient $  \varepsilon$ w.r.t. the point $  y$. The result is the blue point $  B$.

Bella\"{\i}che claims that when $  \varepsilon$ goes to $  0$ the point $  B$ converges to the sum $  y +_{x} z$. Also, from this intrinsic definition of addition, all the other properties (Carnot group structure) of the operation may be deduced from the uniformity of this convergence. He does not give a proof of this fact.

The idea of Bella\"{\i}che is partially correct (in regards to the emergence of the algebraic properties of the operation from uniformity of the convergence of its definition) and partially wrong (this is not the correct definition of the operation). Let me start with the second part. The definition of the operation has the obvious default that it uses the point $  A$ which is far from $  x$. This is in contradiction with the local character of the definition of the metric tangent space (and in contradiction with the local definition of dilations).  But he is wrong from interesting reasons, as we shall see.

Instead, a slightly different path could be followed, figured by the red points $  C, D, E$. Indeed, instead of going far away first (the blue point $  A$), then coming back at finite distance from $  x$ (the blue point $  B$), we may first come close to $  x$ (by using  the red points $  C, D$), then inflate the point $  D$ to finite distance from $  x$ and get the point $  E$. The recipe is a bit more complicated, it involves three dilations instead of two, but I can prove that it works (and leads to the definition of dilation structures and later to the definition of emergent algebras).

The interesting part is that if we draw, as in the figure here,  the constructions in the euclidean plane, then we get $  E = B$, so actually in this case there is no difference between the outcome of these constructions. At further examination this looks like an affine feature, right? But in fact this is true in non-affine situations, for example in the case of intrinsic dilations in Carnot groups, see the examples from the section \ref{secemergent}.

Let's think again about these dilations, which are central to our discussion, as being operations. We may change the notations like this:

$  \delta^{x}_{\varepsilon} y = x \circ_{\varepsilon} y$

Then, it is easy to verify that the equality between the red point $  E$ and the blue point $  B$ is a consequence of the fact that in usual vector spaces (as well as in their non-commutative version, which are Carnot groups), the dilations, seen as operations, are self-distributive! That is why Bella\"{\i}che is actually right in his definition of the tangent space addition operation, provided that it is used only for self-distributive dilation operations. (But this choice limits the applications of his definition of addition operation only to Carnot groups).

\paragraph{ Remark:} I was sensible to these two last sections of Bella\"{\i}che's paper because I was prepared by one of my previous obsessions, namely how to construct differentiability only from topological data.  This was the subject of my first paper, see the story told in the post "\href{http://chorasimilarity.wordpress.com/2011/06/25/topological-substratum-of-the-derivative/}{Topological substratum of the derivative}", there is still some mystery to it, see \href{http://arxiv.org/abs/0911.4619}{arXiv:0911.4619}.

\section{The groupoid construction}
\label{secthree}

Recall that in section \ref{secone}   is mentioned a paragraph by Bella\"{\i}che, where he explains that the structure of the tangent bundle is "concealed in dilations". This is now clear, I hope, but there is something left to do. Bella\"{\i}che also mentions the classical construction of the tangent bundle made by Connes, where the tangent bundle of the space $  X$  appears as a completion of the trivial pair groupoid $  X \times X$.  What I want to explain now is how the dilations interact with the trivial groupoid to give the tangent bundle, in a sort of generalization of Connes construction. See\href{http://arxiv.org/abs/1107.2823}{ arXiv:1107.2823} for all details.

\paragraph{ 1. The trivial pair groupoid} over a set $  X$ is the set $  X \times X$ with  the partially defined operation:

$  (x,u) (u,v) = (x,v)$

This partially defined operation is the composition of arrows in the groupoid which has $  X$ as the set of objects and $  X \times X$ as the set of arrows.  The source of the arrow $  (x,y)$ is $  y = \alpha(x,y)$, the target of that arrow is $  x = \omega(x,y)$.  (This may seem strange, but we have to define the source and target like this in order to  see $  (x,u) (u,v) = (x,v)$ as composition of arrows in a groupoid).

The inverse of the arrow $  (x,y)$ is $  (y,x) = (x,y)^{-1}$.

The addition of arrows (i.e. composition) is

$  add [(x,u), (u,v)] = (x,v)$

and it is defined for pairs of arrows $  ((x,u), (u,v))$ such that $\alpha(x,u) = \omega(u,v)$.

It will be useful further to introduce the difference of two arrows:

$  dif [(u,x), (v,x)] = (u,x) \, (v,x)^{-1} = (u,x) (x,v) = (u,v)$

The difference is a partially defined operation which  makes sense (is defined for) for pairs of arrows with the same source.

\paragraph{ 2. Dilations from the groupoid point of view. }Let $  \Gamma = (0,+\infty)$ be the multiplicative group of strictly positive reals.  For any $  \varepsilon \in \Gamma$ we define the {\em dilation of arrows} of coefficient $  \varepsilon$ to be:

$  \delta_{\varepsilon} (x,y) = (\delta_{\varepsilon}^{y} x , y)$

where $  \delta_{\varepsilon}^{y} x$ is the intrinsic dilation of coefficient $  \varepsilon$, of the point $  x$ with respect to the basepoint $  y$.  Again, as in the case of the definition of source and target of an arrow, here the definition of the dilation of arrows is turned on its head with respect to the writing convention from left to right.

For simplicity we don't care about the domain of definition of the intrinsic dilations, or about the domain of definition of the dilation of arrows. Enough is to say that the dilation of arrows is defined for small enough arrows. What could that mean ?It is simple: think about the distance function $  d: X \times X \rightarrow [0,+\infty)$ as if it is defined on the trivial groupoid. Then, the distance function appears as a kind of a norm on the set of arrows of the groupoid. The norm of the arrow $  (x,y)$ is simply $  d(x,y)$  (first time remarked by Lawvere). Then, a short enough arrow is one with a small norm.

Dilations of arrows have the following properties:
\begin{enumerate}
	\item[-]they preserve the source of arrows:  $  \alpha (\delta_{\varepsilon} (x,y)) = y = \alpha(x,y)$, 
	\item[-]upon the identification $  x \equiv (x,x)$ of objects with their respective identity arrows, any dilation of arrows preserves the objects: $  \delta_{\varepsilon}(x,x) = (x,x)$, as a consequence of the fact that $  \delta_{\varepsilon}^{x}x = x$, 
	\item[-]finally, they form a one-parameter group: $  \delta_{\varepsilon} ( \delta_{\mu} (x,y)) = \delta_{\varepsilon \mu}(x,y)$, because of the relation $  \delta^{x}_{\varepsilon} \delta^{x}_{\mu} y = \delta^{x}_{\varepsilon \mu} y$. 
\end{enumerate}

\paragraph{ 3. Deforming the trivial groupoid. }  I want to use the trivial groupoid in order to explain the recipe for the approximate sum, given by the red dots in the figure from section \ref{sectwo}. 

 With the  notations from emergent algebras explained in section \ref{secemergent} we see that $  E = \Sigma^{x}_{\varepsilon} (y,z)$.  There is an equivalent construction for the approximate difference of two points, with respect to a basepoint:

\vspace{.5cm}   \centerline{\includegraphics[width=120mm]{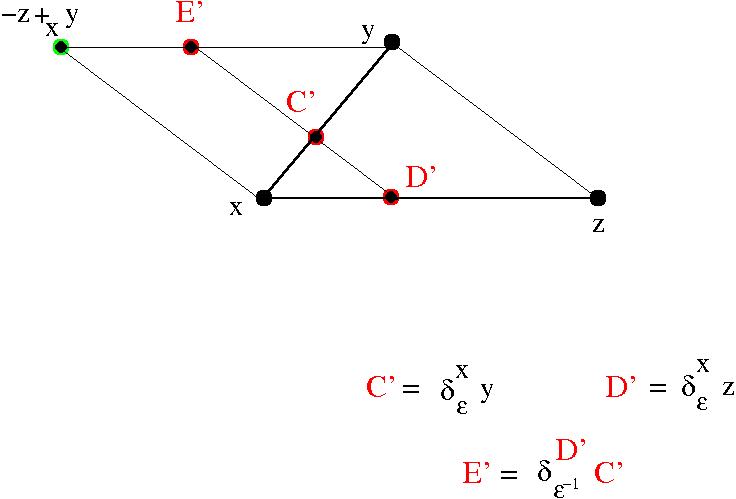}}  \vspace{.5cm}

Again with the notations from emergent algebras we see that $  E' = \Delta^{x}_{\varepsilon}(z,y)$.

By using the trivial groupoid and its deformations by dilations of arrows, it  is easier to explain the origin of the approximate difference. That is because the difference of two arrows is defined on pairs of arrows which have the same source and because dilations of arrows preserve the sources of arrows.

Take two arrows with the same source, say $  g = (y,x)$ and $  h = (z,x)$. Then $  g_{\varepsilon} = \delta_{\varepsilon}(y,x)$ and $  h_{\varepsilon} = \delta_{\varepsilon}(z,x)$ have the same source. We can define then the difference of those deformed arrows:

$$  dif(g_{\varepsilon}, h_{\varepsilon}) = (\delta^{x}_{\varepsilon} y , \delta^{x}_{\varepsilon} z)$$

We define now the deformed difference $  dif_{\varepsilon}$, which is defined on pairs of arrows with the same source, by the relation

$$  \delta_{\varepsilon} dif_{\varepsilon} (g,h) = dif(g_{\varepsilon}, h_{\varepsilon}) $$

A short computation gives

$$  dif_{\varepsilon} [ (z,x) , (y,x) ] = \, ( \Delta^{x}_{\varepsilon}(z,y) , \delta_{\varepsilon}^{x} z) $$

or, with the notations from the last figure

$$  dif_{\varepsilon} [ (z,x) , (y,x) ] = \, (E', D')$$

Finally we can pass to the limit with $  \varepsilon$, as in the tangent bundle construction by Connes.

\section{Emergent algebras}
\label{secemergent}

\begin{definition}  Let $  \Gamma$  be a commutative group with neutral element denoted by $  1$ and operation denoted multiplicatively. A $  \Gamma$ idempotent quasigroup is a set $  X$ endowed with a family of operations $  \circ_{\varepsilon}: X \times X \rightarrow X$,  indexed by $  \varepsilon \in \Gamma$, such that:
\begin{enumerate}
	\item[-]For any $  \varepsilon \in \Gamma$ the pair $  (X, \circ_{\varepsilon})$ is an \href{http://en.wikipedia.org/wiki/Idempotent}{idempotent} \href{http://en.wikipedia.org/wiki/Quasigroup#Definitions}{quasigroup}, 
	\item[-]The operation $  \circ_{1}$ is trivial: for any $  x,y \in X$ we have $  x \circ_{1} y = y$, 
	\item[-]For any $  x, y \in X$ and any $  \varepsilon, \mu \in \Gamma$ we have:  $  x \circ_{\varepsilon} ( x \circ_{\mu} y) = x \circ_{\varepsilon \mu} y$. 
\end{enumerate}
\label{defi1}
\end{definition}

This definition may look strange, let me give some examples of $  \Gamma$ idempotent quasigroups.

\paragraph{ Example 1.}  Real vector spaces: let $  X$ be a real vector space, $  \Gamma = (0,+\infty)$ with multiplication of reals as operation. We define, for any $  \varepsilon > 0$ the following quasigroup operation:

$$  x \circ_{\varepsilon} y = (1-\varepsilon) x + \varepsilon y$$

These operations give to $  X$ the structure of a $(0,+\infty)$ idempotent quasigroup.  Notice that $  x \circ_{\varepsilon}y $ is the dilation based at $  x$, of coefficient $  \varepsilon$, applied to $  y$.

\paragraph{  Example 2.} Complex vector spaces: if $  X$ is a complex vector space then we may take $  \Gamma = \mathbb{C}^{*}$ and we continue as previously, obtaining an example of a $  \mathbb{C}^{*}$ idempotent quasigroup.

\paragraph{  Example 3.} Contractible groups: let $  G$ be a group endowed with a group morphism $  \phi: G \rightarrow G$. Let $  \Gamma = \mathbb{Z}$ with the operation of addition of integers (thus we may adapt Definition 1 to this example by using "$  \varepsilon + \mu$" instead of "$  \varepsilon \mu$" and "$  0$" instead of "$  1$" as the name of the neutral element of $  \Gamma$).  For any $  \varepsilon \in \mathbb{Z}$ let

$$  x \circ_{\varepsilon} y = x \phi^{\varepsilon}(x^{-1} y)$$

This a $  \mathbb{Z}$ idempotent quasigroup. The most interesting case (relevant also for Definition 3 below) is the one when $  \phi$ is an uniformly contractive automorphism of the topological group $  G$. The structure of these groups is an active exploration area, see for example \href{http://arxiv.org/abs/0704.3737}{arXiv:0704.3737}by  \href{http://arxiv.org/find/math/1/au:+Glockner_H/0/1/0/all/0/1}{Helge Glockner}   and the bibliography therein  (a fundamental result here is Siebert article Contractive automorphisms on locally compact groups, Mathematische Zeitschrift 1986, Volume 191, Issue 1, pp 73-90).  See also conical groups and relations between contractive and conical groups introduced in \href{http://arxiv.org/abs/0804.0135}{arXiv:0804.0135},  shortly explained in \href{http://arxiv.org/abs/1005.5031}{arXiv:1005.5031}.

\paragraph{  Example 4.} \href{http://en.wikipedia.org/wiki/Carnot_group}{Carnot groups}: these are a particular example of a conical group. The most trivial noncommutative Carnot group is the \href{http://en.wikipedia.org/wiki/Heisenberg_group}{Heisenberg group}.

\paragraph{ Example 5.} A group with an invertible self-mapping $  \phi: G \rightarrow G$  such that $  \phi(e) =e$, where $  e$ is the identity of the group $  G$. In this case the construction from Example 3 works here as well because there is no need for $  \phi$ to be a group morphism.

\paragraph{  Example 6.} Local versions. We may accept that there is a way (definitely needing care to well formulate, but intuitively cleart) to define a local version of the notion of a $  \Gamma$  idempotent quasigroup. With such a definition, for example, a convex subset of a real vector space gives a local $(0,+\infty)$ idempotent quasigroup (as in Example 1) and a neighbourhood of the identity of a topological group $  G$, with an identity preserving, locally defined invertible self map (as in Example 5) gives a $  \mathbb{Z}$ local idempotent quasigroup.

\paragraph{  Example 7.} A particular case of Example 6, is a Lie group $  G$ with the operations  defined for any $  \varepsilon \in (0,+\infty)$ by

$$  x \circ_{\varepsilon} y = x \exp ( \varepsilon \log (x^{-1} y) )$$

\paragraph{  Example 8. } A less symmetric example is the one of $  X$ being a riemannian manifold, with associated operations  defined for any $  \varepsilon \in (0,+\infty)$ by

$$  x \circ_{\varepsilon}y = \exp_{x}( \varepsilon \log_{x}(y))$$

\paragraph{  Example 9. } More generally, any metric space with dilations  (introduced in  \href{http://arxiv.org/abs/math/0608536}{arXiv:math/0608536}[MG] )  is a local idempotent quasigroup.

\paragraph{ Example 10. }  One parameter deformations of quandles. A \href{http://en.wikipedia.org/wiki/Racks_and_quandles}{quandle} is a self-distributive quasigroup. Take now a one-parameter family of quandles (indexed by $  \varepsilon \in \Gamma$) which satisfies moreover points 2. and 3. from Definition 1. What is interesting about this example is that quandles appear as decorations of knot diagrams, which are preserved by the Reidemeister moves.  At closer examination, examples 1-4 are all particular cases of one parameter quandle deformations!

\vspace{.5cm}

I shall define now the operations of approximate sum and approximate difference associated to a $  \Gamma$  idempotent quasigroup.

For any $  \varepsilon \in \Gamma$, let use define $  x \bullet_{\varepsilon} y = x \circ_{\varepsilon^{-1}} y$.

\begin{definition}  For any $  \varepsilon \in \Gamma$ we give the following names to several combinations of operations of emergent algebras:
\begin{enumerate}
	\item[-]the approximate sum operation is $  \Sigma^{x}_{\varepsilon} (u,v) = $ $  x \bullet_{\varepsilon} ((x \circ_{\varepsilon} u) \circ_{\varepsilon} v)$, 
	\item[-]the approximate difference operation is $  \Delta^{x}_{\varepsilon} (u,v) = (x \circ_{\varepsilon} u) \bullet_{\varepsilon} (x \circ_{\varepsilon} v)$ 
	\item[-]the approximate inverse operation is $  inv^{x}_{\varepsilon} u = (x \circ_{\varepsilon} u) \bullet_{\varepsilon} x$. 
\end{enumerate}
\label{defi2}
\end{definition}

Suppose now that $  X$ is a separable uniform space.  

Let us suppose that the commutative group $\Gamma$ is a topological group endowed with an absolute, i.e. with an invariant topological filter, denoted by $  0$. We write $  \varepsilon \rightarrow 0$ for a net in $ \Gamma$ which converges to the filter $  0$. The image to have in mind is $  \Gamma = (0, + \infty)$ with multiplication of reals as operation and with the filter $  0$ as the filter generated by sets $  (0, a)$ with $  a > 0$. This filter is the restriction to the set $  (0,\infty) \subset \Gamma$ of the filter of  neighbourhoods of the number $  0 \in \mathbb{R}$.  Another example is $  \Gamma = \mathbb{Z}$ with addition of integers as operation, seen as a discrete topological group, with the absolute generated by sets $  \left\{ n \in \mathbb{Z} \, \, : \, \, n \leq M \right\}$ for all $  M \in \mathbb{Z}$. For this example the neutral element (denoted by $  1$ in Definition 1) is the integer $  0$, therefore in this case we can change notations from multiplication to addition, $  1$ becomes $  0$, the absolute $  0$ becomes $  - \infty$ , and so on.

\begin{definition}  An emergent algebra (or uniform idempotent quasigroup) is a $  \Gamma$ idempotent quasigroup  $  X$, as in Definition \ref{defi1}, which satisfies the following topological conditions:
\begin{enumerate}
	\item[-]The family of operations $  \circ_{\varepsilon}$ is compactly contractive, i.e. for any compact set $  K \subset X$, for any $  x \in K$ and for any open neighbourhood $  U$ of $  x$, there is an open set $  A(K,U) \subset \Gamma$ which belongs to the absolute $  0$ such that for any $  u \in K$ and $  \varepsilon \in A(K,U)$ we have $  x \circ_{\varepsilon} u \in U$. 
	\item[-]As $  \varepsilon \rightarrow 0$ there exist the limits 
\end{enumerate}
$  \lim_{\varepsilon \rightarrow 0} \Sigma^{x}_{\varepsilon} (y,z) = \Sigma^{x} (y,z)$  and $  \lim_{\varepsilon \rightarrow 0} \Delta^{x}_{\varepsilon} (y,z) = \Delta^{x} (y,z)$

and moreover these limits are uniform with respect to $  x,y,z$ in compact sets.
\label{defi3}
\end{definition}

The justification for these names comes from the explanations given in the section \ref{secthree}, where I discussed the sketch of a solution to the question "What makes the (metric)  tangent space (to a sub-riemannian regular manifold) a group?", given by Bella\"{\i}che in the last two sections of his article  The tangent space in sub-riemannian geometry, in the book Sub-riemannian geometry, eds. A. Bella\"{\i}che, J.-J. Risler, Progress in Mathematics 144, Birkhauser 1996. We have seen there that the group operation (the noncommutative,  in principle, addition of vectors) can be seen as the limit of compositions of intrinsic dilations, as $  \varepsilon$ goes to $  0$. It is important that this limit exists and that it is uniform, according to Gromov's hint.

The structure theorem of emergent algebras is the following:

\begin{theorem} Let $  X$ be  a $  \Gamma$ emergent algebra. Then for any $  x \in X$ the pair $  (X, \Sigma^{x}(\cdot, \cdot))$ is a conical group.
\label{theo1}
\end{theorem}

Well,  with the notation $  \delta^{x}_{\varepsilon} y = x \circ_{\varepsilon} y$, $  \delta^{x}_{\varepsilon^{-1}} y = x \bullet_{\varepsilon} y$, it becomes clear, for example, that the composition of intrinsic dilations described in the figure from the section  \ref{secthree} is nothing but the approximate sum from Definition \ref{defi2}. (This is to say that formally, if we replace the emergent algebra operations with the respective intrinsic dilations, then the approximate sum operation $  \Sigma^{x}_{\varepsilon}(y,z)$  appears as the red point E from the mentioned  figure. It is still left to prove that intrinsic dilations from regular sub-riemannian spaces give rise to emergent algebras, this was done in \href{http://arxiv.org/abs/0708.4298}{arXiv:0708.4298}.)

We recognize therefore the two ingredients of Bella\"{\i}che's solution into the definition of an emergent algebra:
\begin{enumerate}
	\item[-]approximate operations, which are just clever compositions of intrinsic dilations in the realm of sub-riemannian spaces, which 
	\item[-]converge in a uniform way to the exact operations which give the algebraic structure of the tangent space. 
\end{enumerate}
Therefore, a rigorous formulation of Bella\"{\i}che's solution is Theorem 1 from the previous post, provided that we extract,  from the long differential geometric work done by Bella\"{\i}che, only the part which is necessary for proving that intrinsic dilations produce an emergent algebra structure.

Nevertheless, Theorem \ref{theo1} shows that the "emergence of operations" phenomenon is not at all specific to sub-riemannian geometry. In fact, once we get the idea of the right definition of approximate operations (from sub-riemannian geometry), we can simply try to prove the theorem by "abstract nonsense", i.e. algebraically, with a dash of uniform convergence at the end.

For this we have to identify the algebraic relations which are satisfied by these approximate operations.  For example, is the approximate sum associative? is the approximate difference the inverse of the approximate sum? is the approximate inverse of an element the inverse with respect to the approximate sum? and so on. The answer to these questions is "approximately yes".

It is clear that in order to find the right relations (approximate associativity and so on) between these approximate operations we need to reason in a more clear way. Just by looking at the expressions of the operations from Definition \ref{defi2}, it is obvious that if we start with a brute force  "shut up and compute" approach  then we will end rather quickly with a mess of parantheses and coefficients. There has to be a more easy way to deal with those approximate operations than brute force.

The way I have found has to do with a graphical representation of these operations, a way which eventually led me to \href{http://chorasimilarity.wordpress.com/graphic-lambda-calculus/}{graphic lambda calculus}.

\end{document}